\documentclass[journal]{IEEEtran}

\usepackage{epsf}
\usepackage{epsfig}
\usepackage{cite}
\usepackage{graphicx}
\usepackage{url}
\usepackage{amsmath}

\usepackage{color}

\begin{document}
\title{Experimental evidence of variable-order behavior \\
of ladders and nested ladders
\thanks
{
Manuscript was prepared on July 11, 2011. 
This work was partially supported by the Polish Ministry of
Science and Higher Education grant number 4125/B/T02/2009/36
and the European Union in the framework of European Social
Fund through the Warsaw University of Technology Development
Programme (by Center for Advanced Studies WUT).
This work was also supported in parts by the Slovak Grant Agency for Science
under grants VEGA 1/0390/10, 1/0497/11, 1/0746/11, and also by
grants APVV-0040-07, SK-FR-0037-09, SK-PL-0052-09, and SK-AT-0024-10.
}
}
\author{Dominik Sierociuk, ~\IEEEmembership{Member, IEEE}, Igor Podlubny, ~\IEEEmembership{Member, IEEE}, Ivo~Petr\'a\v{s},~\IEEEmembership{Member, IEEE}
\thanks{Dominik Sierociuk is with Institute of Control and Industrial Electronics, Warsaw University of Technology, Koszykowa 75, Warsaw, Poland,
        E-mail: dsieroci@isep.pw.edu.pl.}
\thanks{
        Igor Podlubny and Ivo Petr\'a\v{s} are with
        Institute of Control and Informatization of Production Processes,
        BERG Faculty, Technical University of Ko\v{s}ice,
        B. N\v{e}mcovej 3, 042 00 Ko\v{s}ice, Slovakia,
        Tel.: +421-55-602-(5179, 5194), E-mail: igor.podlubny@tuke.sk, ivo.petras@tuke.sk}
}
\date{}

\maketitle

\begin{abstract}
The experimental study of two kinds of electrical circuits, a domino ladder and 
a nested ladder, is presented. While the domino ladder is known and already appeared
in the theory of fractional-order systems, the nested ladder circuit is presented
in this article for the first time. 

For fitting the measured data, a new approach is suggested, which is based on using the 
Mittag-Leffler function and which means that the data are fitted by a solution of an
initial-value problem for a two-term fractional differential equation. 

The experiment showed that  in the frequency domain the domino ladder
behaves as a system of order 0.5 and the nested ladder as a system of order 0.25, 
which is in perfect agreement with the theory developed for their design. 

In the time domain, however, the order of the domino ladder is changing 
from roughly 0.5 to almost 1, and the order of the nested ladder is changing
in a similar manner, from roughly 0.25 to almost 1; in both cases, the 
order 1 is never reached, and both systems remain the systems of non-integer order less than 1.

Both studied types of electrical circuits provide the first known
examples of circuits, which are made of passive elements only
and which  exhibit in the time domain the behavior of variable order.
\end{abstract}

\begin{IEEEkeywords}
fractional calculus, variable order, fractance, fractional integrator,
domino ladder, nested ladder,  Mittag-Leffler function
\end{IEEEkeywords}

\maketitle

\section{Introduction to Fractional Calculus}
\label{sec1}

\PARstart{F}{ractional} calculus is more then 300 years old topic, 
which during recent decades became a powerful and widely used tool
for better modeling and control of processes in many fields of science
and engineering \cite{Podlubny,Oldham,Magin-book,blas2010,Caponetto}.
The term ``fractional calculus" has some historical background and is used
for denoting the theory of integration and differentiation of 
arbitrary real (not necessarily integer) order.  

The standard notation for denoting the left-sided fractional-order differentiation 
of a function $f(t)$ defined in the interval $[a, b]$ is $_{a}D^{\alpha}_{t} f(t)$,
 with $\alpha \in \mbox{R}$.  
 Sometimes a simplified notation $f^{(\alpha)}(t)$ or $d^\alpha f(t)/dt^\alpha$ is used.
 In some applications also right-sided
 fractional derivatives $_{t}D^{\alpha}_{b} f(t)$ are used, but in the present article
 we will use only left-sided fractional derivatives.  Even from the notation one can see 
that evaluation of the left-sided fractional-order operators require the values 
of the function $f(t)$ in the interval $[a, t]$. When $\alpha$ becomes an integer number,
this interval shrinks to the vicinity of the point $t$,  and we obtain the classical
integer-order derivatives as particular cases.

There are several definitions of the fractional derivatives and integrals,
of which we need only the following two.

The Caputo definition of fractional differentiation can be written as \cite{Podlubny}:
\begin{equation}\label{Caputo}
   _{a}^CD_{t}^{\alpha}f(t)=
    \frac{1}{\Gamma (n-\alpha)}
    \int_{a}^{t}
    \frac{f^{(n)}(\tau)}{(t-\tau)^{\alpha - n + 1}}d\tau,
\end{equation}
$$
(n-1 \leq \alpha <n)
$$
\noindent
 where $\Gamma(z)$ is Euler's gamma function.

Above Caputo definition is extremely useful in the time domain studies, because
the initial conditions for the fractional-order differential equations with the Caputo derivatives can be given in the same form as for the integer-order differential equations.
This is an advantage in applied problems, which require the use of initial conditions
containing starting values of the function and its integer-order derivatives
$f(a)$, $f^{'}(a)$, $f^{''}(a)$, \dots, $f^{(n-1)}(a)$.

The formula for the Laplace transform of the Caputo fractional
derivative (\ref{Caputo}) has the form \cite{Podlubny}:
\begin{equation}\label{LTFC}
      \int_{0}^{\infty} e^{-st}\, _{0}^CD_{t}^{\alpha}f(t) \, dt =
      s^{\alpha}F(s) - \sum_{k=0}^{n-1} s^{\alpha-k-1} \,
                     f^{(k)}(0),
\end{equation}
$$
(n-1 \leq \alpha <n).
$$

If the process $f(t)$ is considered from the state of absolute rest, 
so  $f(t)$ and its integer-order derivatives are equal to zero at  
the starting time $t=0$, then the Laplace transform of the $\alpha$-th
derivative of $f(t)$ is simply $s^{\alpha}F(s)$.

The second definition, which we need, is the definition of 
the left-sided Caputo-Weyl fractional derivative:
 \begin{equation}\label{Weyl}
   _{-\infty}^{~~W}D_{t}^{\alpha}f(t)=
    \frac{1}
           {\Gamma (n-\alpha)}
    \int_{-\infty}^{t}
    \frac{f^{(n)}(\tau)d\tau}
    {(t-\tau)^{\alpha -n+ 1} },
\end{equation}
$$
(n-1 \leq \alpha <n)
$$

The Fourier transform of $_{-\infty}^{~~W}D_{t}^{\alpha}f(t)$   
is simply  $(j\omega)^\alpha$.
The Caputo-Weyl definition must be used in the frequency domain studies of 
fractional-order systems. The Caputo-Weyl derivative can be considered
as the Caputo derivative with \hbox{$a \rightarrow -\infty$}.
In other words, the Caputo definition allows the study of the transient
effects in fractional-order systems, which were initially at the state of rest,
while the Caputo-Weyl definition allows the study of frequency responses
of such systems.

Fractional-order models have been already used for modeling 
of electrical circuits (such as domino ladders, tree structures, etc.) 
and elements (coils, memristor, etc.). The review of such models 
can be found in \cite{petras2002,Petras2010mem,Petras2011}. 

Let us consider, for instance, 
a capacitor as a basic element of many circuits.
Westerlund and Ekstam in 1994 proposed a new linear capacitor model
\cite{Westerlund2}. It is based on Curie's empirical law of 1889
which states that the current through a capacitor is
$$
i(t)=\frac{u_0 }{h_1 t^{\alpha}},
$$
where $h_{1}$ and $\alpha$ are constant, $u_{0}$ is the \textit{dc}
voltage applied at $t = 0$, and $0 < \alpha < 1$, \,$(\alpha \in
\mbox{R})$.

For a general input voltage $u(t)$ the current is
\begin{equation}
\label{eq5} i(t)=C \frac{d^{\alpha} u(t)}{dt^{\alpha}} \equiv C \,
_0D^{\alpha}_{t} u(t),
\end{equation}
where $C$ is capacitance of the capacitor. It is related to the kind
of dielectric used in the capacitor. 
The order $\alpha$ is related to the
losses of the capacitor. Westerlund and Ekstam provided in their work the table
of various capacitor dielectrics with appropriated constants $\alpha$
which have been obtained experimentally by measurements.

The relationship between the current and the voltage 
in a~capacitor is described using fractional-order integration:  
\begin{equation}
u(t)=\frac{1}{C} \int_0^t i(t) dt^{\alpha} \equiv \frac{1}{C}\,
_0D^{-\alpha}_{t} i(t).
\end{equation}

Then the impedance of a fractional capacitor is:
\begin{equation}\label{Zc}
Z_c(s) = \frac{1}{C\, s^{\alpha}} = \frac{1}{\omega^{\alpha} C} e^{j(-\alpha\frac{\pi}{2})}, \quad
\omega \in (-\infty, \infty).
\end{equation}

\begin{figure}[h]
\centerline{\includegraphics[width=0.4\textwidth]{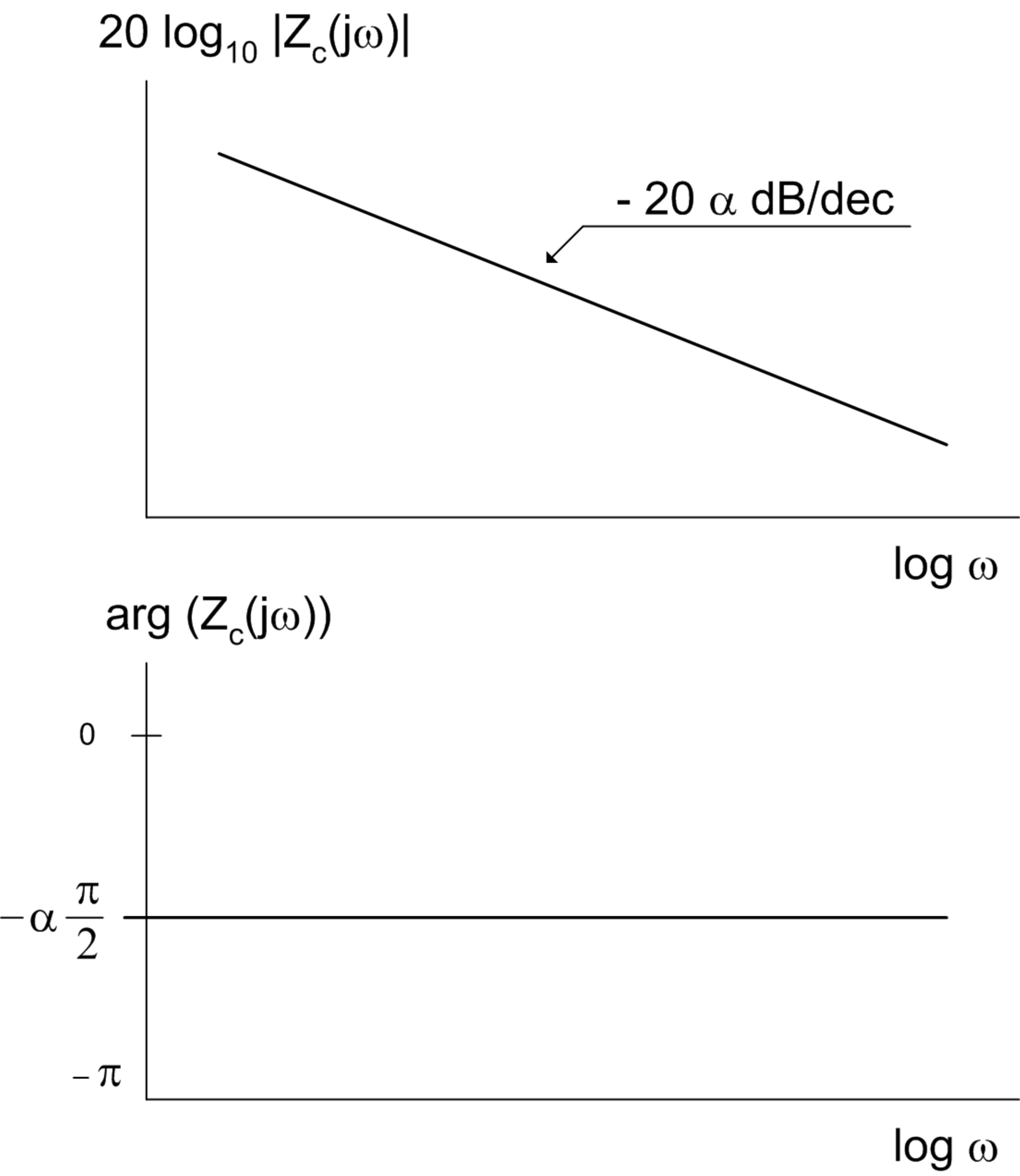}}
\caption{Bode plots of real capacitor.}
\label{BZc}
\end{figure}

Ideal Bode's characteristics of the transfer function for a~real capacitor (\ref{Zc}) are
depicted in Fig.~\ref{BZc}.

General characteristics of the transfer function of a real capacitor
 (\ref{Zc}) are \cite{PetrasETFA2009}:
\begin{itemize}
\item  Magnitude: constant slope of $-\alpha20dB/dec.;$
\item  Crossover frequency: a function of $1/C$;
\item  Phase: horizontal line of $-\alpha \frac{\pi }{2};$
\end{itemize}

Besides this fractional-order capacitor model, 
we can mention the new fractional-order models of coils \cite{Schafer-k2},
memristive systems \cite{Coopmans-Petras-Chen-k2}, ultracapacitors \cite{Dzelinski, Dzelinski2}, 
and the element called fractor \cite{Bohannan}. Such elements can be
combined with classical passive and active elements for creating
various types of electrical circuits.

Among the aforementioned fractional-order elements, the fractor is of special 
interest, because it is known that the order of fractor slowly changes in time
with aging of chemical materials of which it is composed \cite[Table I]{Sheng2011}. In other words,
fractor is an example of an element of variable non-integer order. 
Such variable-order behavior of the fractor was experimentally studied in \cite{Sheng2011}.
In this paper we
demonstrate that variable-order behavior can be observed
in a wide class of  ladder-type circuits composed of standard passive elements. 

\section{Fractional Devices and Fractance}
\label{sec2}

Besides simple elements like a capacitor, electrical circuits of more or less
complex structure were studied by many authors. The review of most of the previous
efforts can be found in \cite{petras2002}. 
A circuit that exhibits fractional-order behavior is called a \emph{fractance}~\cite{Podlubny}. 

\subsection{Fractances}

The fractance devices have the following characteristics \cite{Nakagawa}. 
First, the phase angle is constant independent of the frequency 
within a wide frequency band. 
Second, it is possible to construct a filter having a moderate characteristics 
which can not be realized by using the conventional devices.

Generally speaking, there are three basic types of fractances.
The most popular is a domino ladder circuit network \cite{roy1967}. 
Another type is a tree structure of electrical elements \cite{Nakagawa}, 
and finally, we can find out also some transmission line circuit 
(or symmetrical domino ladder \cite{Carlson2}).

Design of fractances having given order $\alpha$ can be done easily using any of
the  rational approximations or a truncated continued fraction expansion (CFE),
which also gives a rational approximation \cite{Petras-Analog, erfani2002}.
Truncated CFE does not require any further transformation;
a rational approximation based on any other methods must be
first transformed to the form of a continued fraction; then the values
of the electrical elements, which are necessary for building
a fractance, are determined from the obtained finite
continued fraction. If all coefficients of the obtained
finite continued fraction are positive, then the fractance
can be made of classical passive elements (resistors and capacitors).
If some of the coefficients are negative, then the fractance
can be made with the help of negative impedance 
converters~\cite{petras2002,Petras-Analog}.

It is worth mentioning also the constant phase element (CPE), which exhibits
the  fractional-order behavior as well. It is a metal-insulator-solution or metal-insulator-liquid  interface used in electrochemistry. CPE interprets a dipole layer capacitance  \cite{Biswas}. The impedance of CPE is expressed as $Z_{CPE}(s)=Q s^{-\alpha}$ and CPE cannot be described by a finite number of passive elements with frequency independent values.

\subsection{Traditional domino ladder (half-order integrator)}

Several different algorithms for approximation the fractional order integrators 
are currently available \cite{petras2002, roy1967, roy1974, Krishna2008, Yifei2005, Wang}.
Most of them are based on some form of approximation of irrational transfer functions
in the complex domain. The commonly used approaches include the aforementioned CFE method and its modifications, or representation by a quotient of polynomials in $s$ in a factorized form.
  
The main disadvantage of these algorithms is that the values of electrical elements (like resistors and capacitors) 
needed for the approximation are not the standard values of elements
produced by manufacturers. 

However, it is still possible to obtain highly accurate and practically usable 
implementations of a fractional-order integrator using only standard elements 
with the standard values available in the market. The idea of this practical approach 
to implementation of fractional-order systems is based on 
the domino ladder structure.

\begin{figure}[thpb]
\setlength{\unitlength}{1mm}
\begin{center}
\begin{picture}(80,30)(5,0)
\put(8,20){\vector(1,0){8}}
\put(12,23){\makebox(0,0){$i(t)$}}
\put(16,4){\vector(0,1){15}}
\put(13,10){\makebox(0,0){$u(t)$}}
\multiput(10,20)(18,0){3}{\line(10,0){10}}
\multiput(20,22)(18,0){3}{\line(10,0){8}}
\multiput(20,18)(18,0){3}{\line(10,0){8}}
\multiput(20,22)(18,0){3}{\line(0,-1){4}}
\multiput(28,22)(18,0){3}{\line(0,-1){4}}
\multiput(33,20)(18,0){3}{\line(0,-1){8}}
\multiput(30,12)(18,0){3}{\line(1,0){6}}
\multiput(30,11)(18,0){3}{\line(1,0){6}}
\multiput(33,11)(18,0){3}{\line(0,-1){8}}
\multiput(31.5,3)(18,0){3}{\line(1,0){3}}
\put(64,20){\line(10,0){8}}
\put(75,19.8){\makebox(0,0){$\dots$}}
\put(28,12){\makebox(0,0){C}}
\put(46,12){\makebox(0,0){C}}
\put(64,12){\makebox(0,0){C}}
\put(24,24){\makebox(0,0){$R$}}
\put(42,24){\makebox(0,0){$R$}}
\put(60,24){\makebox(0,0){$R$}}
\end{picture}
\end{center}
\caption{Domino ladder scheme.}
\label{fig:domino}
\end{figure}
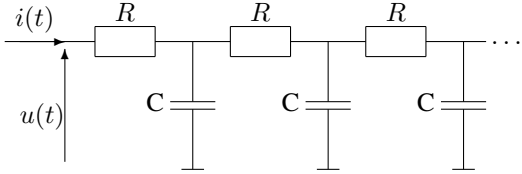

The domino ladder circuit shown in Fig. \ref{fig:domino} has the following impedance:
\begin{equation}\label{DL-standard-impedance}
G(s)=R+\frac{1}{sC +\frac{1}{R+\frac{1}{sC+ \frac{1}{R+ \frac{1}{sC +\frac{1}{R +\frac{1}{sC +\dots }}}}}}} = \frac{1}{(Ts)^{0.5}},
\end{equation}
where $T={C/R}$. In the ideal case of infinite realization, (\ref{DL-standard-impedance})
gives a half-order integrator; a truncated realization gives its approximation.

The domino ladder circuit can be also considered as a model of a semi-infinite RC line, which is described by the following partial differential equations \cite{Manabe, Oldham73}:
\begin{equation}
\frac{\partial}{\partial x}u(t,x)=R i(t,x),\label{eqn:par_1}
\end{equation}
\begin{equation}
\frac{\partial}{\partial x}i(t,x)=C\frac{\partial}{\partial t}u(t,x),\label{eqn:par_2}
\end{equation}
where $u(t,x)$ is the voltage and $i(t,x)$ is the current at point $x$ at time instance $t$.

This can be rewritten as 
\begin{equation}
\frac{\partial^2}{\partial x^2}u(t,x)=RC\frac{\partial}{\partial t}u(t,x).
\end{equation}

From this equation a relationship between  the current $i(t,0)$ 
and voltage $u(t,0)$ at the beginning of the semi-infinite RC line 
can be obtained in terms of half-order integral; 
in the Laplace domain it has the following form:

\begin{equation}
G(s)=\frac{U(s,0)}{I(s,0)}= \sqrt{\frac{R}{C}} \frac{1}{ s^{0.5}},  
\end{equation}
where $I(s,0)$ and $U(s,0)$ are the Laplace transforms of $i(t,0)$ and $u(t,0)$.

\subsection{Domino ladder with alternating resistors}
\label{sec3}

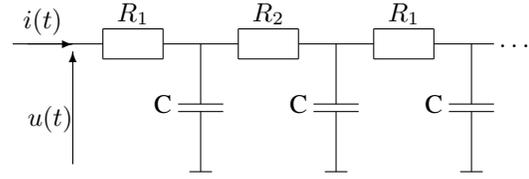
\begin{figure}[thpb]
\setlength{\unitlength}{1mm}
\begin{center}
\begin{picture}(80,30)(5,0)
\put(8,20){\vector(1,0){8}}
\put(12,23){\makebox(0,0){$i(t)$}}
\put(16,4){\vector(0,1){15}}
\put(13,10){\makebox(0,0){$u(t)$}}
\multiput(10,20)(18,0){3}{\line(10,0){10}}
\multiput(20,22)(18,0){3}{\line(10,0){8}}
\multiput(20,18)(18,0){3}{\line(10,0){8}}
\multiput(20,22)(18,0){3}{\line(0,-1){4}}
\multiput(28,22)(18,0){3}{\line(0,-1){4}}
\multiput(33,20)(18,0){3}{\line(0,-1){8}}
\multiput(30,12)(18,0){3}{\line(1,0){6}}
\multiput(30,11)(18,0){3}{\line(1,0){6}}
\multiput(33,11)(18,0){3}{\line(0,-1){8}}
\multiput(31.5,3)(18,0){3}{\line(1,0){3}}
\put(64,20){\line(10,0){8}}
\put(75,19.8){\makebox(0,0){$\dots$}}
\put(28,12){\makebox(0,0){C}}
\put(46,12){\makebox(0,0){C}}
\put(64,12){\makebox(0,0){C}}
\put(24,24){\makebox(0,0){$R_1$}}
\put(42,24){\makebox(0,0){$R_2$}}
\put(60,24){\makebox(0,0){$R_1$}}
\end{picture}
\end{center}
\caption{Proposed analogue model of half-order integrator.}
\label{fig:half}
\end{figure}

For building accurate analog approximation of the half-order integrator using easily accessible elements available in the market, 
the approach presented in Fig. \ref{fig:half} can be used.

Based on the observation made in article \cite{SierociukMMAR11}, 
we can formulate the following design algorithm:
\begin{enumerate}
\item[(a)] Choose the values of $R_1$ and $C$ in order to obtain the required 
low frequency limit.
\item[(b)] Choose value of $R_2$ in order to satisfy the condition $R_1\approx 4R_2$. This condition allows to select those values of resistors that are available as manufactured.
\item[(c)] Choose the ladder length $n$ (number of steps in the domino ladder) 
in order to obtain the desired frequency range of approximation.
\end{enumerate}

\subsection{Enhanced domino ladder for half-order integration}

The modified ladder with two alternating values of resistors
performs better than the classical domino ladder, 
but  the phase shift is still equal not to $45^\circ$,
 but to approximately $43^\circ$--$44^\circ$ 
 ($45^\circ$ achieved only at very short frequency range). 

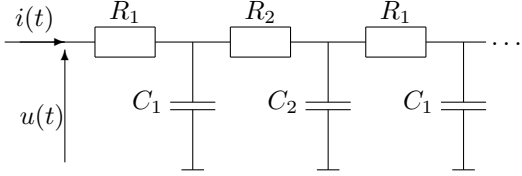
\begin{figure}[thpb]
\setlength{\unitlength}{1mm}
\begin{center}
\begin{picture}(80,30)(5,0)
\put(8,20){\vector(1,0){8}}
\put(12,23){\makebox(0,0){$i(t)$}}
\put(16,4){\vector(0,1){15}}
\put(13,10){\makebox(0,0){$u(t)$}}
\multiput(10,20)(18,0){3}{\line(10,0){10}}
\multiput(20,22)(18,0){3}{\line(10,0){8}}
\multiput(20,18)(18,0){3}{\line(10,0){8}}
\multiput(20,22)(18,0){3}{\line(0,-1){4}}
\multiput(28,22)(18,0){3}{\line(0,-1){4}}
\multiput(33,20)(18,0){3}{\line(0,-1){8}}
\multiput(30,12)(18,0){3}{\line(1,0){6}}
\multiput(30,11)(18,0){3}{\line(1,0){6}}
\multiput(33,11)(18,0){3}{\line(0,-1){8}}
\multiput(31.5,3)(18,0){3}{\line(1,0){3}}
\put(64,20){\line(10,0){8}}
\put(75,19.8){\makebox(0,0){$\dots$}}
\put(27,12){\makebox(0,0){$C_1$}}
\put(45,12){\makebox(0,0){$C_2$}}
\put(63,12){\makebox(0,0){$C_1$}}
\put(24,24){\makebox(0,0){$R_1$}}
\put(42,24){\makebox(0,0){$R_2$}}
\put(60,24){\makebox(0,0){$R_1$}}
\end{picture}
\end{center}
\caption{Enhanced domino ladder for half-order integration.}
\label{fig:quatro}
\end{figure}

To further improve the accuracy of approximation,  let us modify the structure of the domino ladder in such a way that there are not only two values of resistors, 
but also two values of capacitors are used (Fig.~\ref{fig:quatro}).

\begin{figure}[!ht]
\centering
\includegraphics[width=0.48\textwidth]{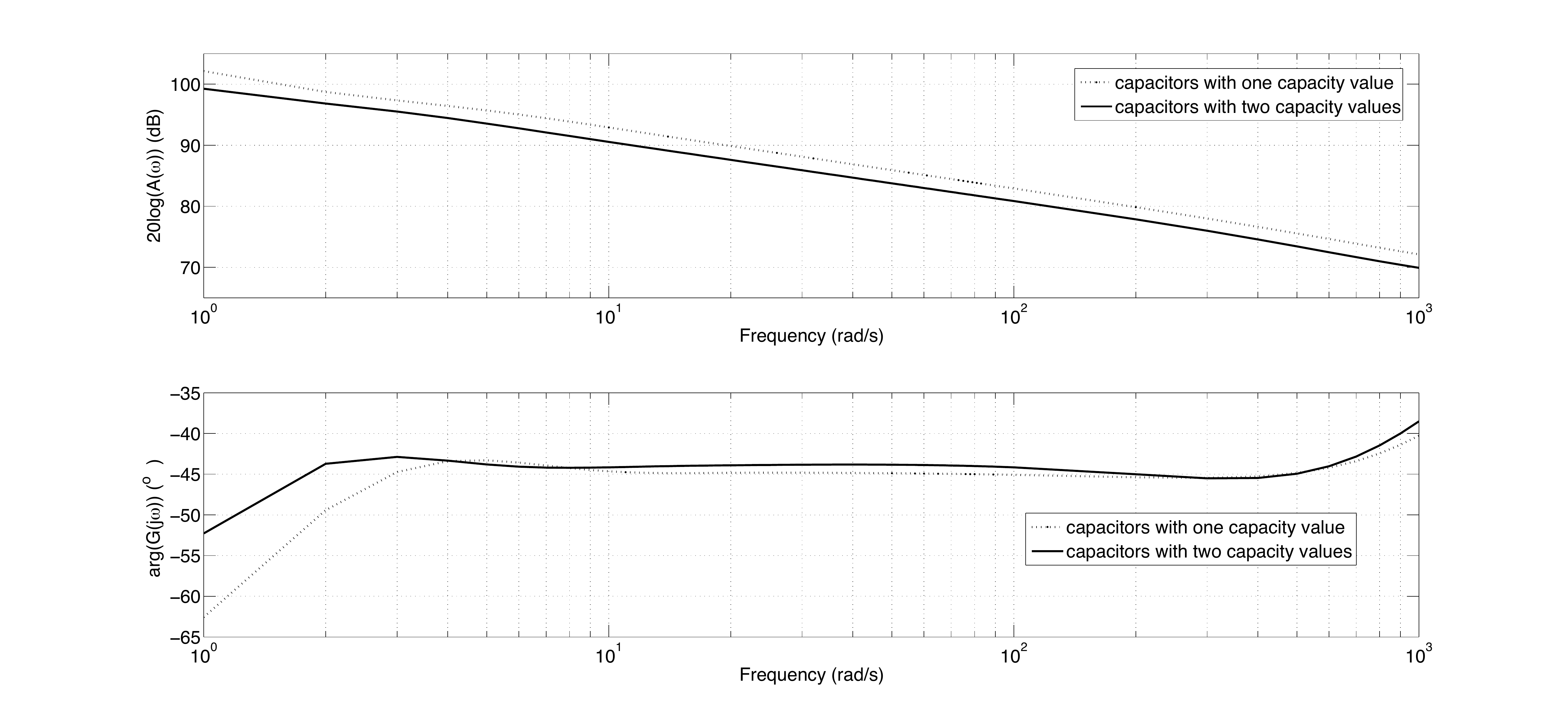}
\caption{Results of modeling of half-order integrator using a modified domino ladder. }
\label{dl_ma}
\end{figure}

Fig \ref{dl_ma} presents the experimental results for the enhanced domino ladder for  half-order integration with the following parameters of the circuit presented in Fig.~\ref{fig:quatro}:
$R_1=2320\Omega$, $R_2=8200\Omega$, $C_1=330{\mathrm{nF}}$, $C_2=220{\mathrm{nF}}$, and the number of steps in the ladder is equal to $n=34$. 
The results are compared with the realization presented in Fig.~\ref{fig:half}.
It is obvious that the phase plot of the enhanced ladder 
is indeed is much closer to the $45^\circ$  than 
in the case of the classical domino ladder.

\subsection{A new type of fractances: a nested ladder}

Based on the above results, we can easily extend them to 
build a fractional order integrator of order $0.25$. 
This can be done by replacing the capacitors in the scheme in Fig.~\ref{fig:half} 
by half-order integrators, which can be either classical domino ladders 
or enhanced domino ladders.
This step can be interpreted as an introduction a half-order dynamics 
into the equation (\ref{eqn:par_2}). 
This results in a transfer function of order $\alpha=0.25$, 
which corresponds to a quarter-order integrator.

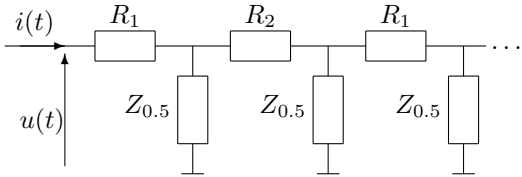
\begin{figure}[thpb]
\setlength{\unitlength}{1mm}
\begin{center}
\begin{picture}(80,30)(5,0)
\put(8,20){\vector(1,0){8}}
\put(12,23){\makebox(0,0){$i(t)$}}
\put(16,4){\vector(0,1){15}}
\put(13,10){\makebox(0,0){$u(t)$}}
\multiput(10,20)(18,0){3}{\line(10,0){10}}
\multiput(20,22)(18,0){3}{\line(10,0){8}}
\multiput(20,18)(18,0){3}{\line(10,0){8}}
\multiput(20,22)(18,0){3}{\line(0,-1){4}}
\multiput(28,22)(18,0){3}{\line(0,-1){4}}
\multiput(33,20)(18,0){3}{\line(0,-1){4}}
\multiput(31,16)(18,0){3}{\line(1,0){4}}
\multiput(31,7)(18,0){3}{\line(1,0){4}}
\multiput(31,16)(18,0){3}{\line(0,-1){9}}
\multiput(35,7)(18,0){3}{\line(0,1){9}}
\multiput(33,7)(18,0){3}{\line(0,-1){4}}
\multiput(31.5,3)(18,0){3}{\line(1,0){3}}
\put(64,20){\line(10,0){8}}
\put(75,19.8){\makebox(0,0){$\dots$}}
\put(27,12){\makebox(0,0){$Z_{0.5}$}}
\put(45,12){\makebox(0,0){$Z_{0.5}$}}
\put(63,12){\makebox(0,0){$Z_{0.5}$}}
\put(24,24){\makebox(0,0){$R_1$}}
\put(42,24){\makebox(0,0){$R_2$}}
\put(60,24){\makebox(0,0){$R_1$}}
\end{picture}
\end{center}
\caption{Integrator of order $\alpha=0.25$ in the form of a nested ladder.}
\label{fig:quadro}
\end{figure}

In Fig.~\ref{fig:quadro} the scheme of the approximation of a quarter-order integrator is shown;  $Z_{0.5}$ are the impedances of modified or enhanced domino ladders implementing half-order integrators. 

In the same way (namely, by replacing impedances $Z_{0.5}$ with $Z_{0.25}$) 
an integrator of order $\alpha=0.125$ and so forth can be built, 
but this will need a large number of elements. 

We call such a structure of electrical circuit the \emph{nested ladder}.
The nested ladder is an example of using the ideas of self-similatiry
and fractality for creating electrical circuits exhibiting non-integer order behavior.

\section{Data fitting using the Mittag-Leffler function}
\label{sec:identification-method}

In order to obtain a model for the measured data from the considered electrical circuits
(ladders and nested ladders),
we have developed a new approach to data fitting, which is based on using the Mittag-Leffler function and which, in fact, allows obtaining models of non-integer order.

The idea of our method is based on the following. 
When it comes to obtaining a mathematical models from measurements or observations,
it is a common practice in many fields of science and engineering to choose the type 
of the fitting curve and identify its parameters using some criterion (usually a least squares method).  We would like to point out that choosing a particular type of a curve means that, in fact, the process is modeled by a differential equation, for which that curve is a solution. 

For example, fitting data using the equation $y(t) = a t + b$ (known as linear regression model) means that the process is modeled by the solution of a simple second-order differential equation under two initial conditions:

\begin{equation}
y'' = 0, \quad
y(0) = b, \quad
y'(0) = a.
\end{equation}

Similarly, the fitting function in the form $y = a \sin(\omega t) + b \cos(\omega t)$ 
means that the process is modeled by the solution of the initial value problem
of the form

\begin{equation}
y'' + \omega^2 y = 0, \quad
y(0) = b, \quad
y'(0) = a \omega.
\end{equation}

Choosing the fitting function in another frequently used form, $y = a e^{bt}$,
means that the process is modeled by the solution of the initial value problem

\begin{equation}
y' - b \, y = 0, \quad
y(0) = a.
\end{equation}

Thinking in this way, we conclude that instead of postulating the shape of the fitting curve
it is possible to postulate the form of the initial-value problem and identify the parameters appearing in the differential equation and in the initial conditions. For the first time this method was suggested and used in \cite[Chapter 10]{Podlubny}. 
In this paper we, however,  just emphasize that obtaining a fitting function $y(t)$
for measurements of a dynamic process immediately means that 
that process is described by an initial-value problem of which $y(t)$ is the solution.

In the present article the measured data are fitted by 
\begin{equation} \label{eq:ML-fitting-function}
y = y_0 \, E_{\alpha, 1} (a\, t^\alpha)
\end{equation}
where $E_{\alpha, \beta}(z)$ is the Mittag-Leffler function defined as~\cite{Podlubny}

\begin{equation}
E_{\alpha, \beta} (z) 
= 
\sum_{k=0}^{\infty}
\frac{z^k}
       {\Gamma (\alpha k + \beta)}.
\end{equation}

The parameters to be identified are $\alpha$,  $a$, and $y_0$.

If the data are fitted by the function (\ref{eq:ML-fitting-function}), then this means 
that they are modeled by the solution of the following initial-value problem
for a two-term fractional-order differential equation containing the Caputo fractional derivative of order~$\alpha$:
\begin{equation}\label{eq:IVP-fitting}
_{0}^{C}\!D_{t}^{\alpha} y(t) - a \, y(t) = 0, \quad y(0) = y_0.
\end{equation}

\section{Experimental Results}
\label{sec5}

\subsection{Experimental setup}

\begin{figure}[!htb]
\centering
\includegraphics[width=0.48\textwidth]{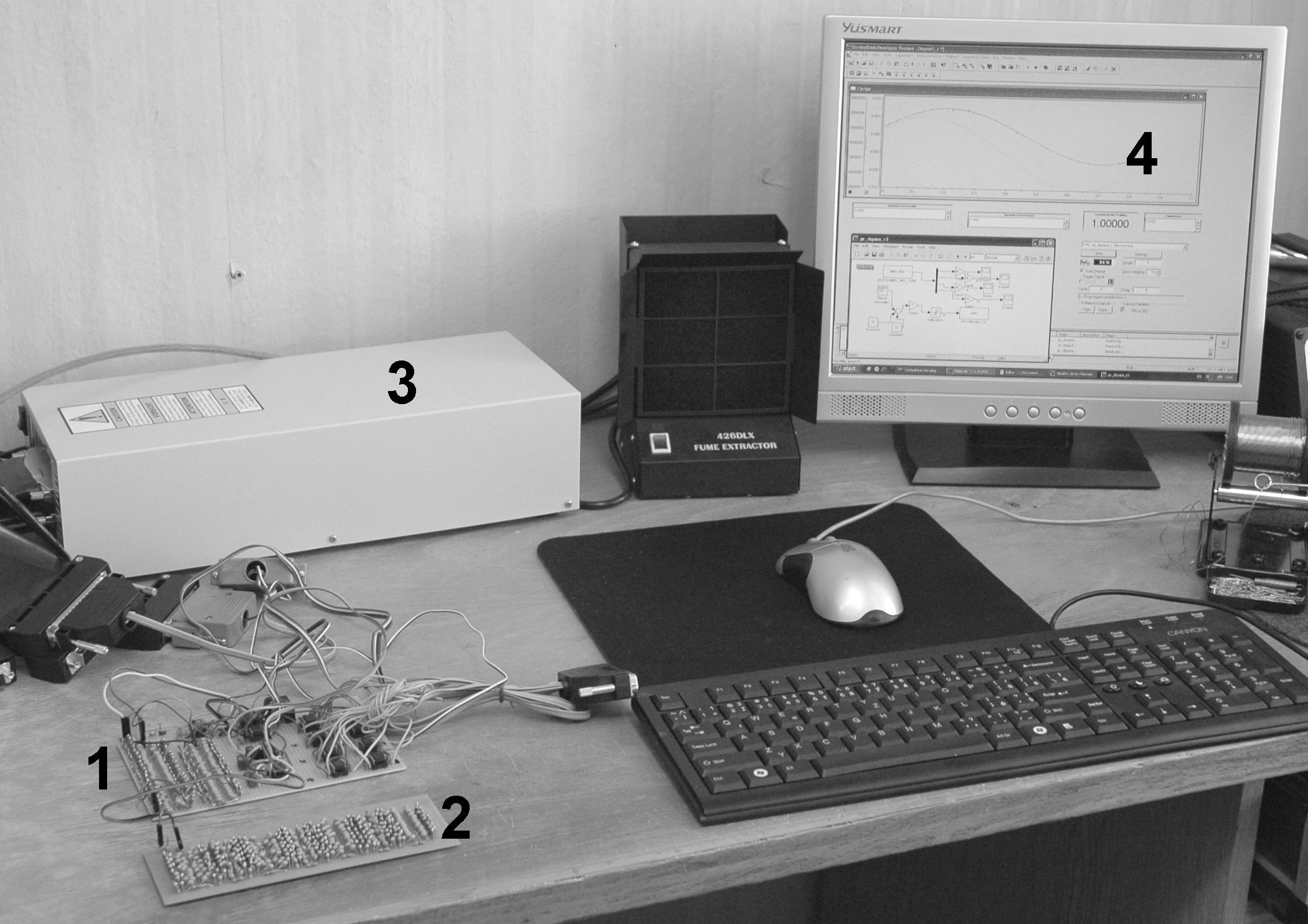}
\caption{Experimental setup used for all measurements: 1 -- domino ladder, 2 -- nested ladder, 3 -- dSpace card, 4 -- computer with Matlab/Simulink software.}
\label{setup}
\end{figure}

\begin{figure}[!htb]
\centering
\includegraphics[width=0.48\textwidth]{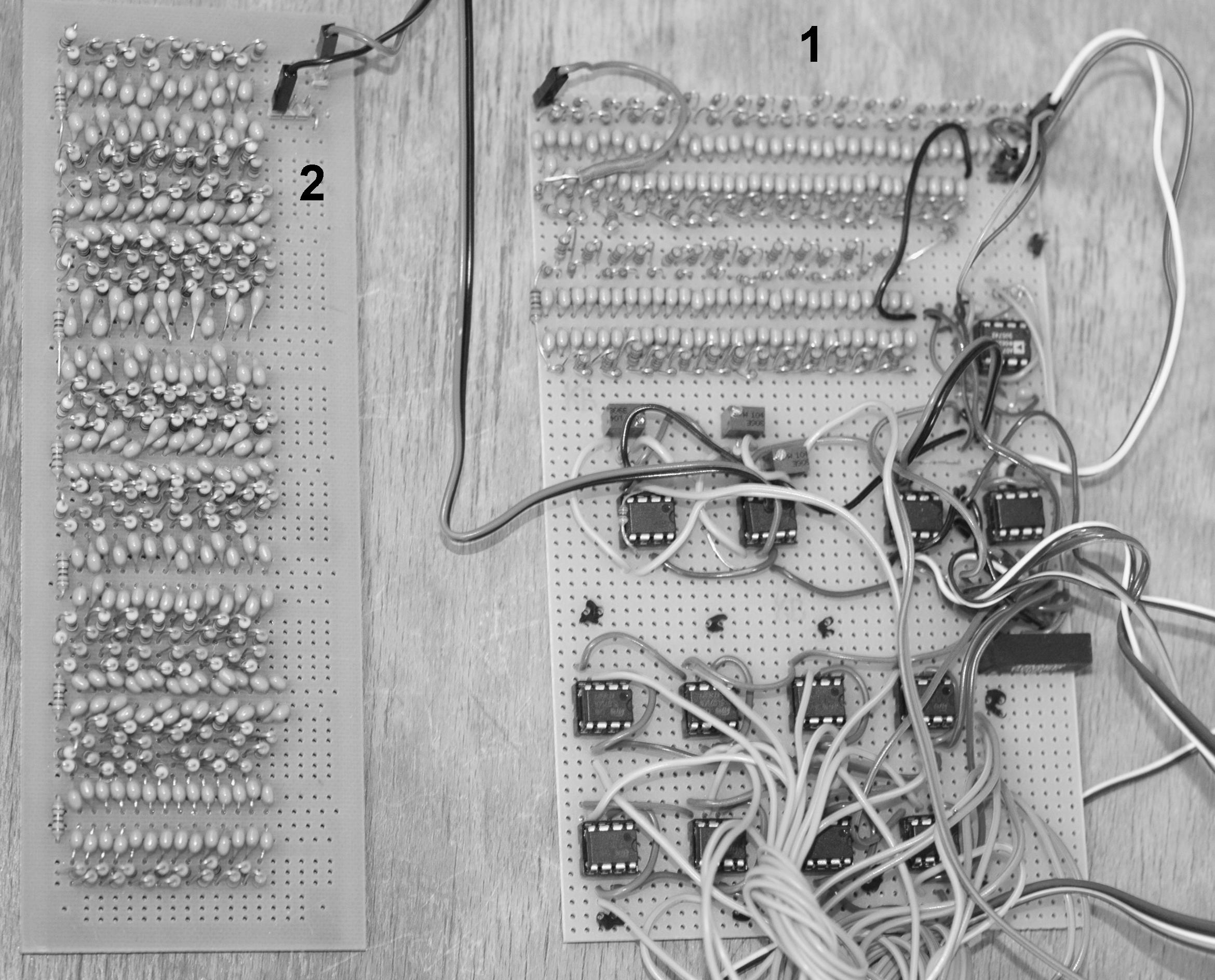}
\caption{Detailed view: 1 -- domino ladder, 2 -- nested ladder. }
\label{setupDetails}
\end{figure}

For the experimental verification of the introduced method,
the circuits presented in Section \ref{sec2}  were built.
For measurements, the modified domino ladder circuit and the nested ladder 
were connected to the amplifier electronic circuit of the operational amplifiers TL071 and to the dSpace DS1103 DSP card connected to a computer. 
The real laboratory setup is shown in Fig.~\ref{setup} and detailed view of the ladders is in Fig.~\ref{setupDetails}.

\begin{figure}[thpb]
\setlength{\unitlength}{0.6mm}
\begin{center}
\begin{picture}(140,40)(-15,0)
\put(15,25){\makebox(0,0){\framebox(10,3){}}}
\put(15,30){\makebox(0,0){{\small{$\mathrm{10k\Omega}$}}}}
\put(27,24){\makebox(0,0){-}}
\put(27,15){\makebox(0,0){+}}
\put(20,25){\line(10,0){5}}
\put(22,25){\line(0,1){8}}
\put(22,33){\line(1,0){8}}
\put(35,33){\makebox(0,0){\framebox(10,3){}}}
\put(38,38){\makebox(0,0){{\small{$\mathrm{Z_{ladder}}$}}}}
\put(40,33){\line(1,0){8}}
\put(48,20){\line(0,1){13}}
\put(3,25){\line(1,0){7}}
\put(2,25){\circle{1}}
\put(2,5){\circle{1}}
\put(2,5){\line(0,-1){3}}
\put(2,8){\vector(0,1){15}}
\put(1,2){\line(1,0){2}}
\put(15,15){\line(10,0){10}}
\put(15,15){\line(0,-1){13}}
\put(13,2){\line(1,0){4}}
\put(25,30){\line(0,-1){20}}
\put(25,30){\line(2,-1){20}}
\put(25,10){\line(2,1){20}}
\put(45,20){\line(10,0){8}}
\put(77,25){\line(0,1){8}}
\put(77,33){\line(1,0){8}}
\put(90,33){\makebox(0,0){\framebox(10,3){}}}
\put(91,38){\makebox(0,0){{\small{${R_i}$}}}}
\put(82,24){\makebox(0,0){-}}
\put(82,15){\makebox(0,0){+}}
\put(95,33){\line(1,0){8}}
\put(103,20){\line(0,1){13}}
\put(70,25){\makebox(0,0){\framebox(10,3){}}}
\put(71,30){\makebox(0,0){{\small{${R_i}$}}}}
\put(75,25){\line(10,0){5}}
\put(70,15){\line(10,0){10}}
\put(70,15){\line(0,-1){13}}
\put(68,2){\line(1,0){4}}
\put(80,30){\line(0,-1){20}}
\put(80,30){\line(2,-1){20}}
\put(80,10){\line(2,1){20}}
\put(100,20){\line(10,0){11}}
\put(53,20){\line(0,1){5}}
\put(53,25){\line(1,0){12}}
\put(-5,15){\makebox(0,0){{\small{${u_1(t)}$}}}}
\put(119,10){\makebox(0,0){{\small{${u_2(t)}$}}}}
\put(112,20){\circle{1}}
\put(112,5){\circle{1}}
\put(112,5){\line(0,-1){3}}
\put(112,8){\vector(0,1){10}}
\put(111,2){\line(1,0){2}}
\end{picture}
\end{center}
\caption{Electronic circuit of measurement setup for integrator.}
\label{schem}
\end{figure}
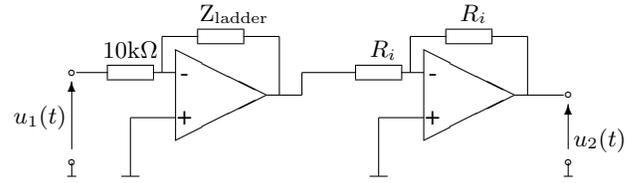

The electronic scheme presented in Fig. \ref{schem} uses two operational amplifiers. The first one is working in the integrator configuration and the second one is working in the inverse unit-gain  for compensating the signal inversion of the integrator amplifier. The resistor $R_i$ can be used for changing the gain of the integrator and it was chosen to $R_i=3.3k\Omega$. The $u_1$ is an input and $u_2$ is an output of the integrator system.

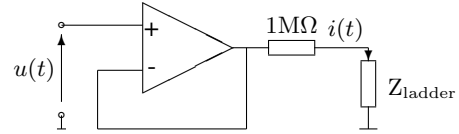
\begin{figure}[thpb]
\setlength{\unitlength}{0.6mm}
\begin{center}
\begin{picture}(80,40)(-2,0)
\put(2,25){\line(10,0){18}}
 \put(10,2){\line(1,0){33}}
 \put(43,2){\line(0,1){18}}
 \put(2,25){\circle{1}}
 \put(2,5){\circle{1}}
 \put(2,5){\line(0,-1){3}}
 \put(2,8){\vector(0,1){15}}
 \put(1,2){\line(1,0){2}}
 \put(10,15){\line(10,0){10}}
 \put(10,15){\line(0,-1){13}}
 \put(20,30){\line(0,-1){20}}
 \put(20,30){\line(2,-1){20}}
 \put(20,10){\line(2,1){20}}
 \put(40,20){\line(10,0){8}}
 \put(53,20){\makebox(0,0){\framebox(10,3){}}}
 \put(58,20){\line(10,0){12}}
 \put(70,20){\vector(0,-1){3}}
 \put(70,7){\line(0,-1){5}}
 \put(70,12){\makebox(0,0){\framebox(3,10){}}}
 \put(68,2){\line(1,0){4}}
\put(82,11){\makebox(0,0){\small$\mathrm{Z_{ladder}}$}}
\put(-4,15){\makebox(0,0){\small{${u(t)}$}}}
\put(22,24){\makebox(0,0){+}}
\put(22,15){\makebox(0,0){-}}
\put(28,20){\makebox(0,0){}}\put(53,25){\makebox(0,0){\small$\mathrm{1M\Omega}$}}
\put(65,24){\makebox(0,0){\small${i(t)}$}}
\end{picture}
\end{center}
\caption{Electronic circuit of direct measurement setup for ladders.}
\label{schemat_dl}
\end{figure}

\subsection{Modified half-order domino ladder measurements}

\begin{figure}[!htb]
\centering
\includegraphics[width=0.48\textwidth]{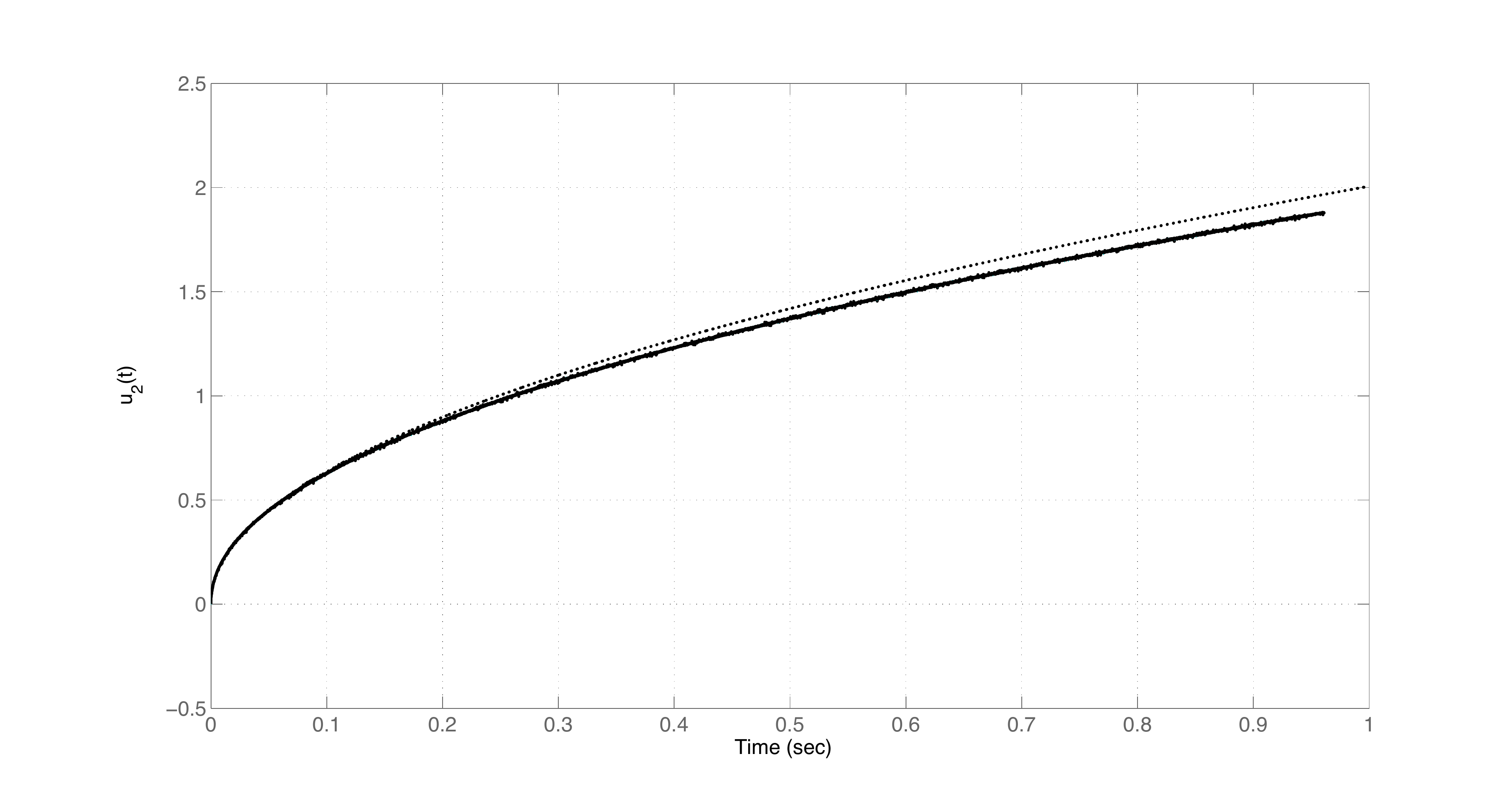}
\caption{Comparison of measured and calculated step responses of half-order integrator with domino ladder of 130 steps:
(dotted line) calculated response for $\alpha=0.5167$ from Table~\ref{tab:VO-table} for $1$ s; (solid line) measured response. }
\label{ex_0_5s}
\end{figure}

\begin{figure}[!ht]
\centering
\includegraphics[width=0.48\textwidth]{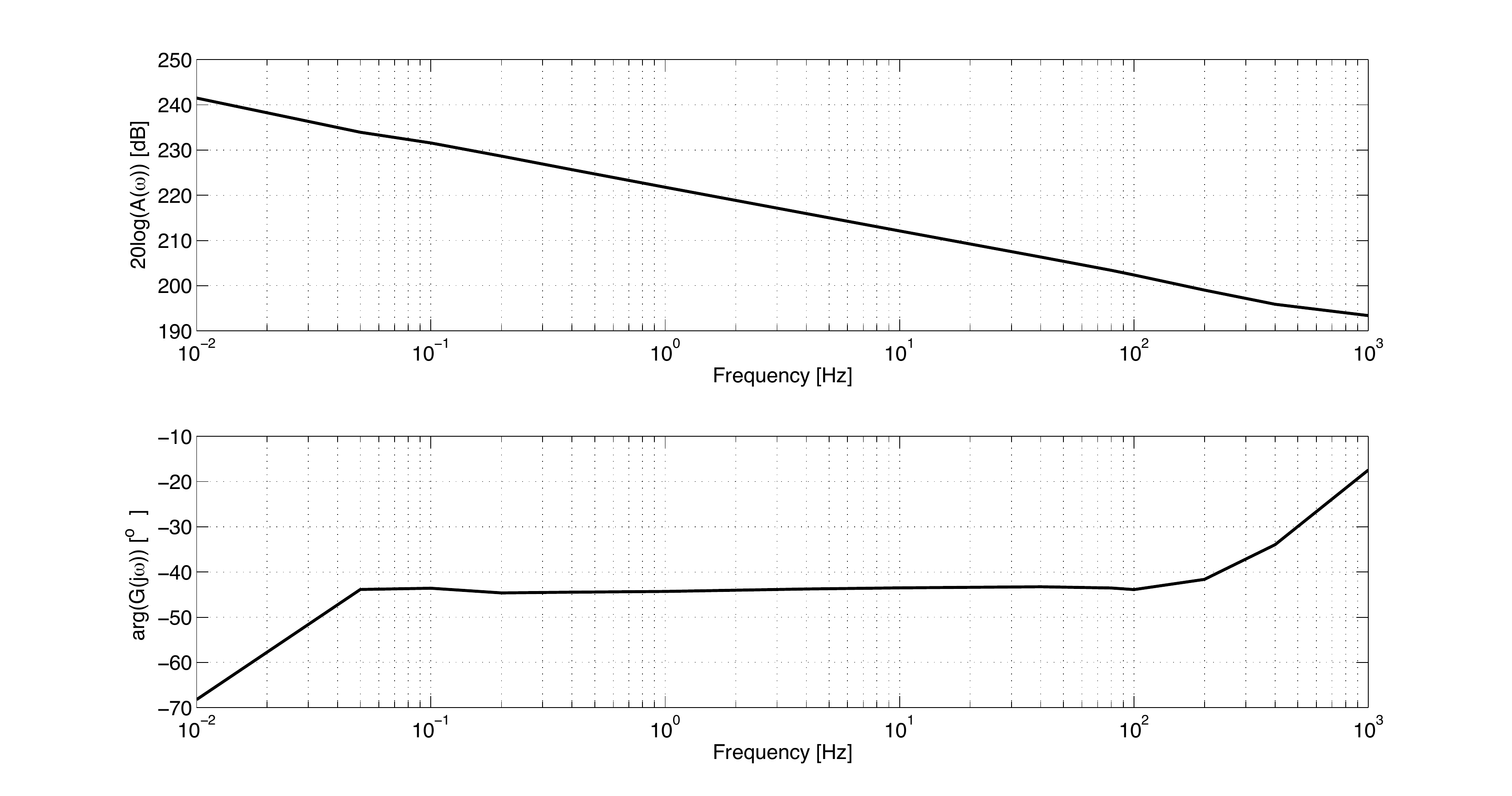}
\caption{Measured Bode plots of half-order integrator with domino ladder of 130 steps. }
\label{ex_0_5b}
\end{figure}

The tested circuit has the following parameters of the circuit presented in Fig. \ref{fig:half}:
$R_1=2000\Omega$, $R_2=8200\Omega$, $C=470{\mathrm{nF}}$, 
and numbers of steps in the ladders was taken first  $n=60$ and then $n=130$.  
The sampling period was $Ts=0.0001$ s.
The manufacturing tolerance of the elements used for making such ladders
is $1\%$ for resistors and $20\%$ for capacitors. 
As it can be seen in Fig.~\ref{ex_0_5s} and Fig.~\ref{ex_0_5b}, the obtained experimental results 
fully confirm the theoretical considerations and simulations. 

\subsection{Quarter-order domino ladder measurements}

\begin{figure}[!htb]
\centering
\includegraphics[width=0.48\textwidth]{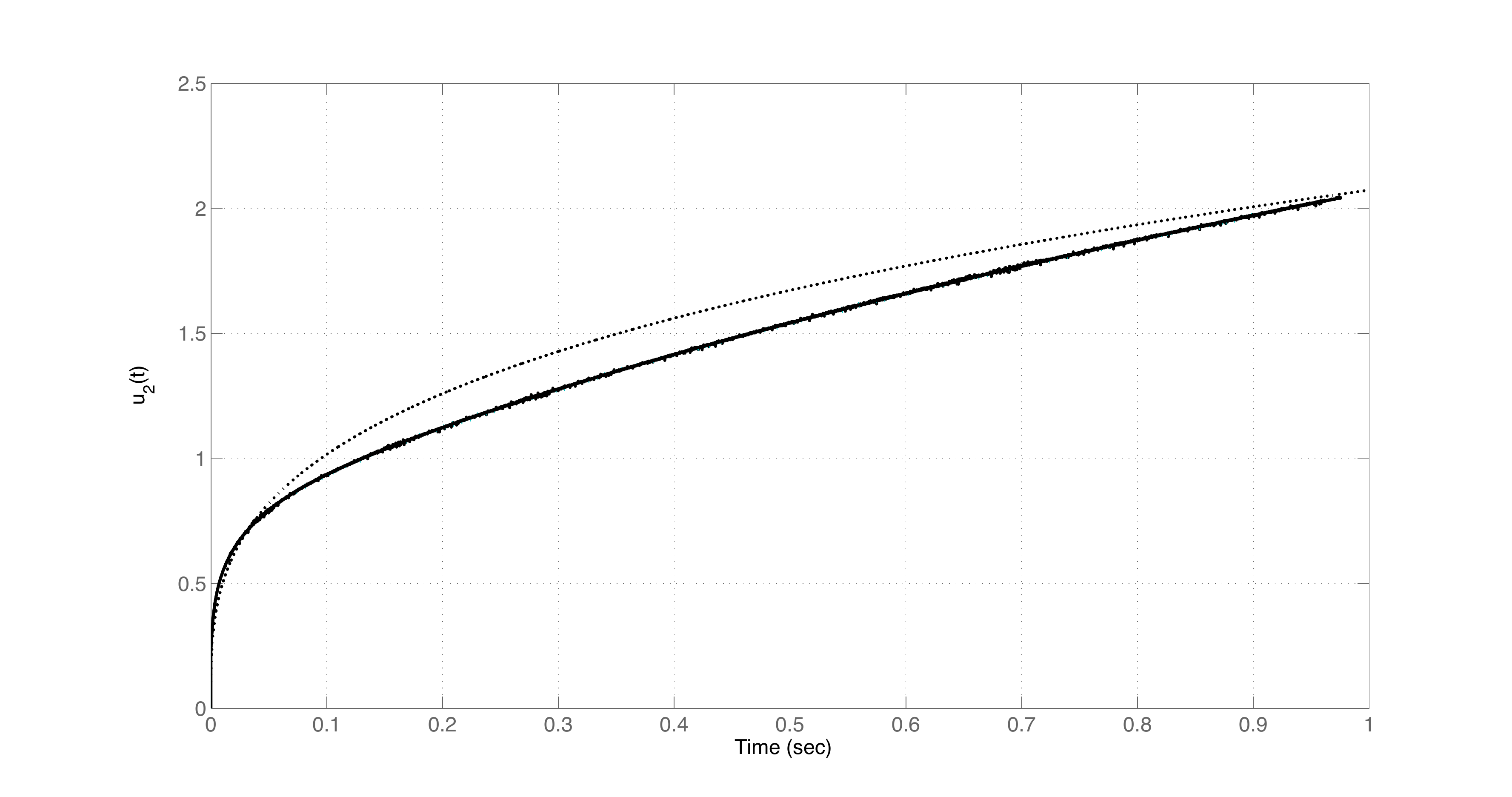}
\caption{Comparison of measured and calculated step responses of quarter-order integrator with nested ladder of size $14 \times 14$ steps: (dotted line) calculated response for $\alpha=0.3126$ from Table~\ref{tab:VO-table} for $1$ s; (solid line) measured response. }
\label{ex_0_25s}
\end{figure}

\begin{figure}[!htb]
\centering
\includegraphics[width=0.48\textwidth]{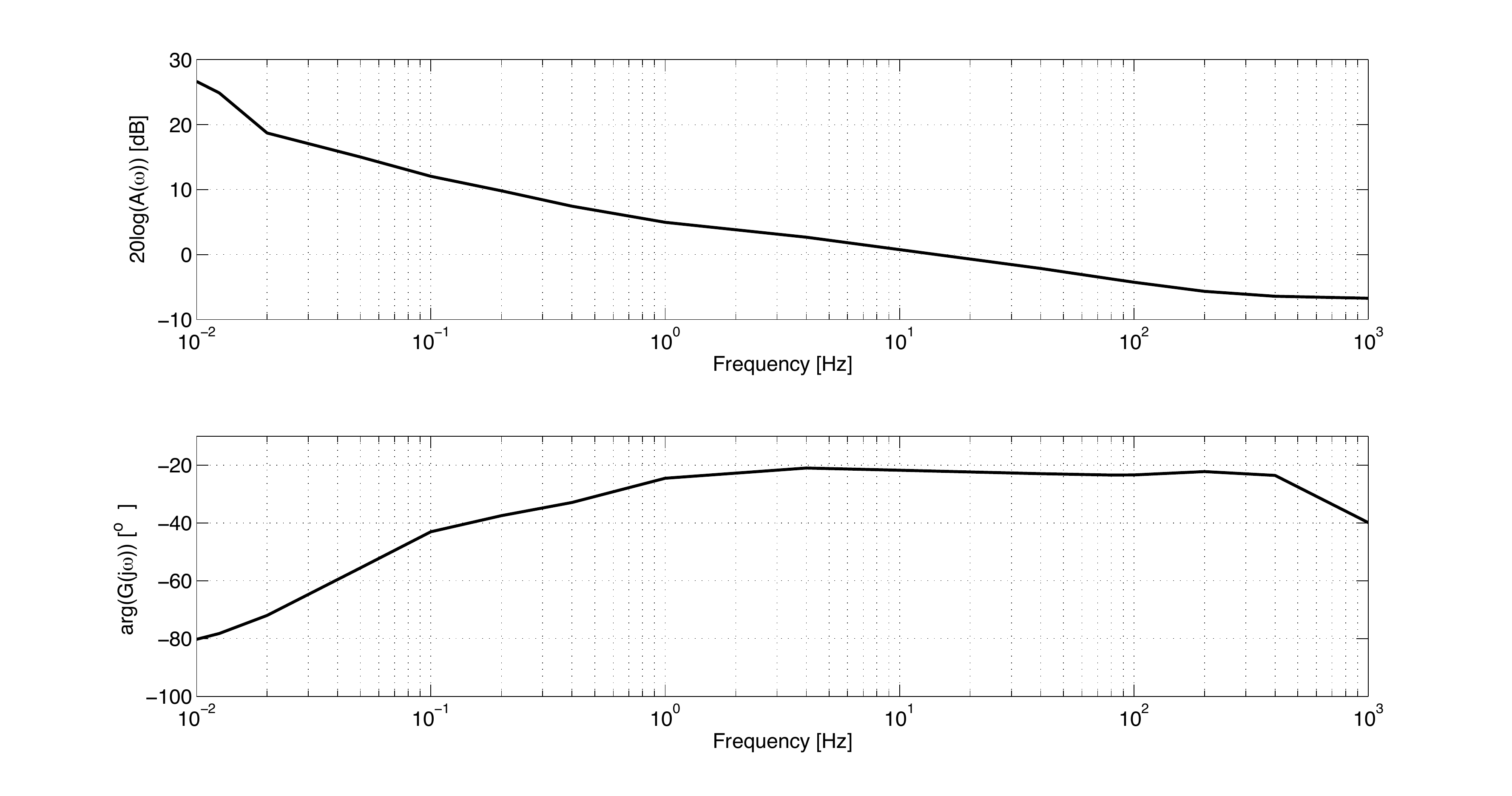}
\caption{Measured Bode plots of quarter-order integrator with nested ladder of size $14 \times 14$ steps. }
\label{ex_0_25b}
\end{figure}

The tested circuit has the following parameters of the circuit presented in Fig. \ref{fig:half}:
$R_1=2000\Omega$, $R_2=8200\Omega$, $C=470{\mathrm{nF}}$ and realization length equal $n=14 \times 14$, that is 14 sub-ladders with 14 steps each.  
The sampling period was $Ts=0.0001$ s.
The manufacturing tolerance of the elements used for making such ladders
is $1\%$ for resistors and $20\%$ for capacitors. 
As it can be seen in Fig.~\ref{ex_0_25s} and Fig.~\ref{ex_0_25b}, the obtained experimental results 
confirm the theoretical considerations and simulations. 
A little deviation in the time domain is due to small number of the nested ladder steps,
which can be also observed in the frequency domain (only two and half decades approximation). 
 
\section{Variable-order behavior}

If the measurements are obtained for the fixed interval $[0, t]$, 
then fitting using the Mittag-Leffler function (\ref{eq:ML-fitting-function}),
described in Section~\ref{sec:identification-method},
immediately gives the model (\ref{eq:IVP-fitting}) of fractional order $\alpha$.

However, if we consider the changing length of the interval, 
then the resulting order of the model will be, in general, a~function of 
this changing interval length $t$: $\alpha = \alpha(t)$.
The same  holds for other two parameters. 

\begin{figure}[htbp]
\begin{center}
\includegraphics[width=\columnwidth]{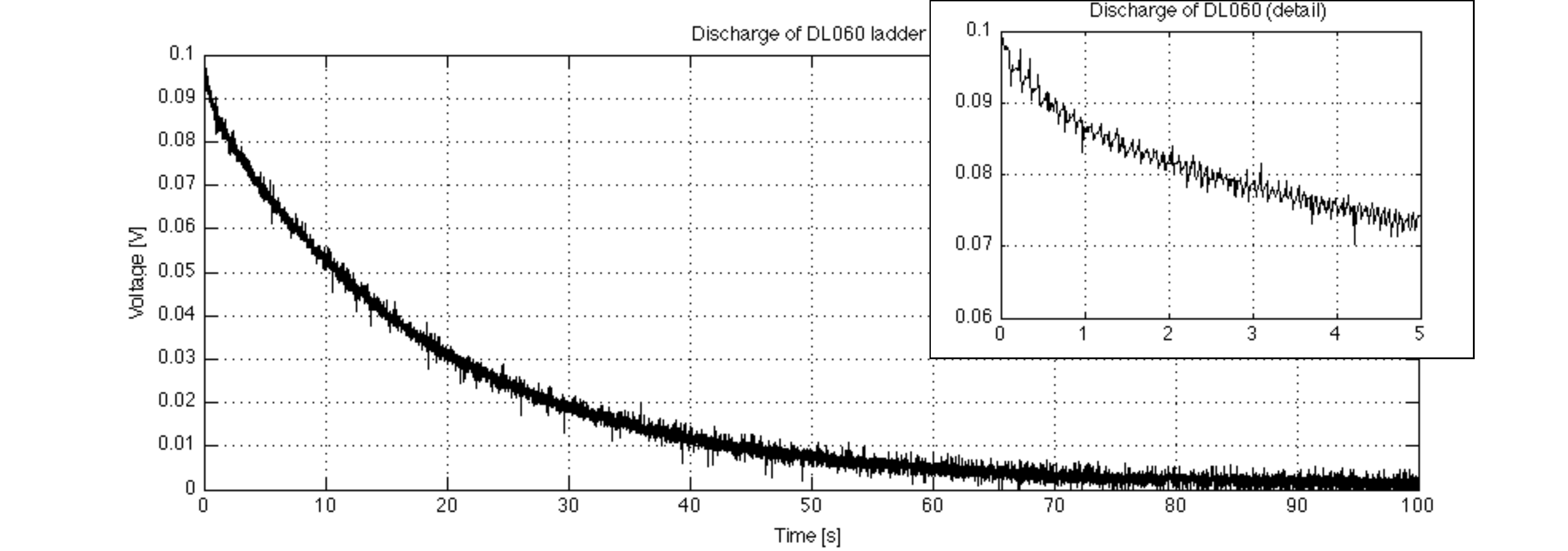}
\caption{Discharge of the  60-steps domino ladder (DL060).}
\label{fig:DL60relax}
\end{center}
\end{figure}

\begin{figure}[htbp]
\begin{center}
\includegraphics[width=\columnwidth]{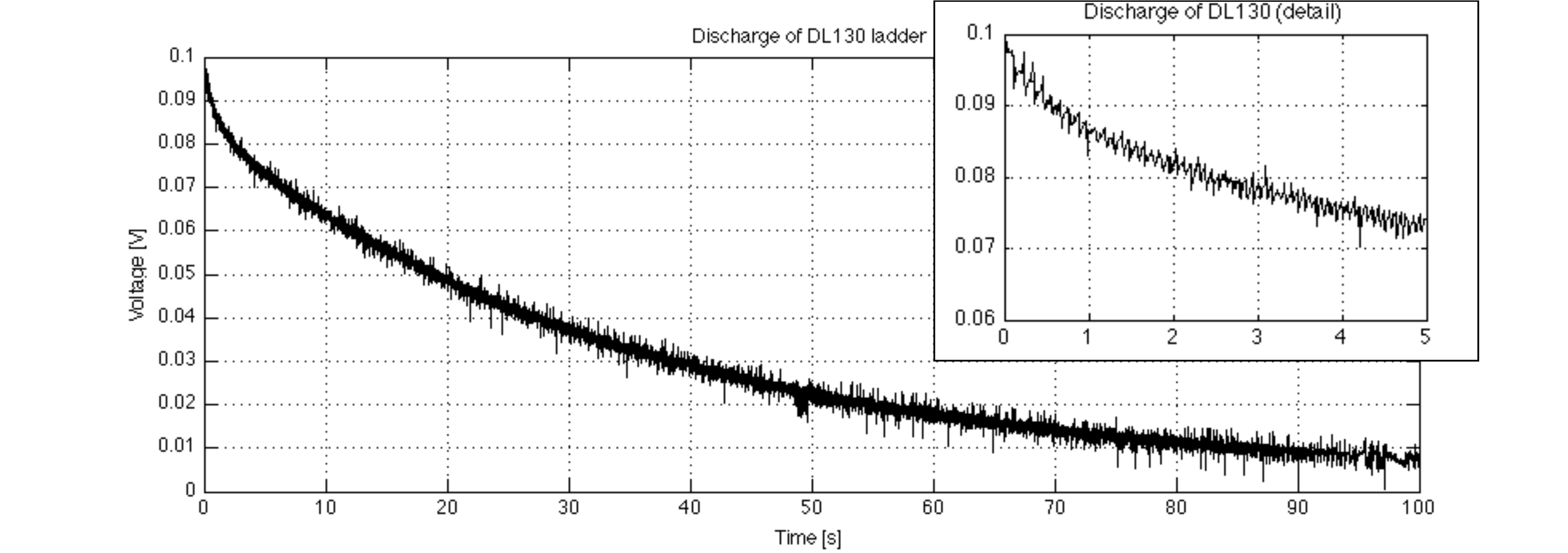}
\caption{Discharge of the 130-steps domino ladder (DL130)).}
\label{fig:DL130relax}
\end{center}
\end{figure}

\begin{figure}[htbp]
\begin{center}
\includegraphics[width=\columnwidth]{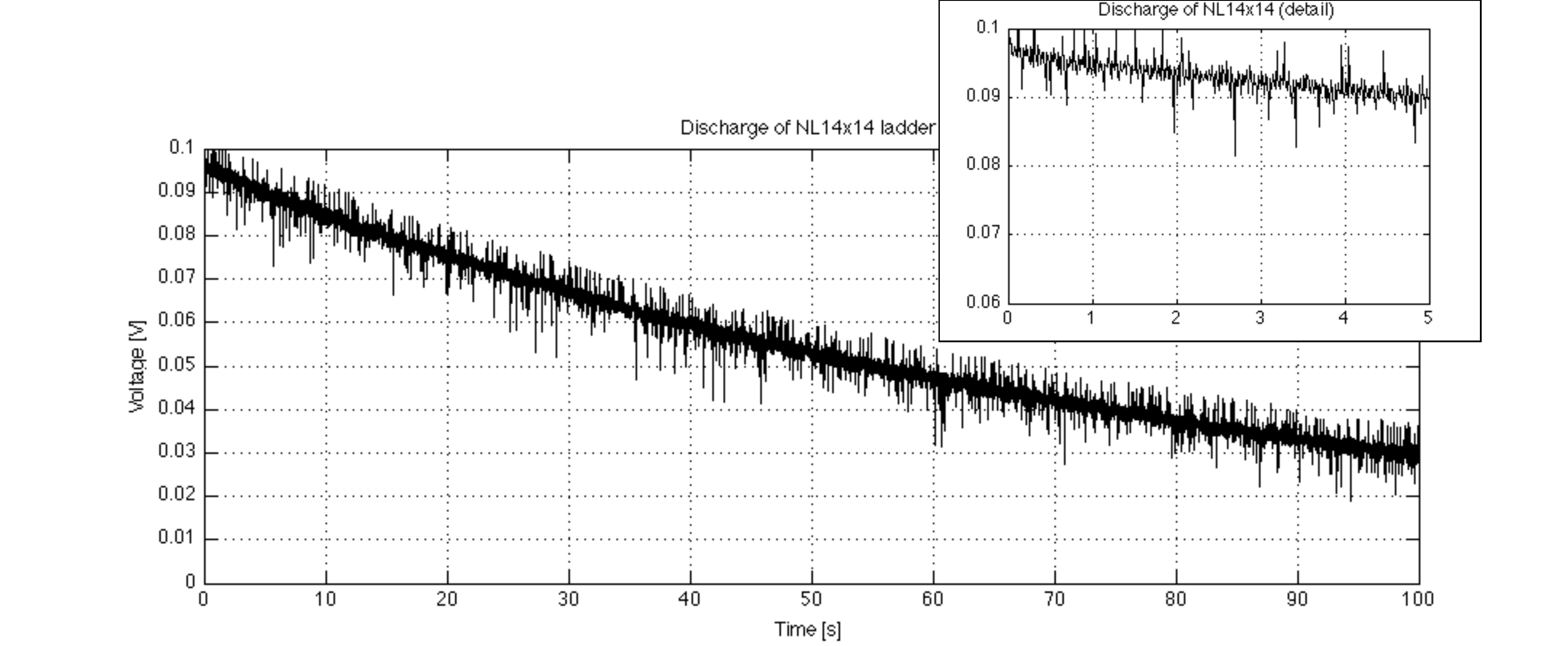}
\caption{Discharge of the nested domino ladder (NL14x14).}
\label{fig:NLrelax}
\end{center}
\end{figure}

In our experiments we considered the growing number of measurements
that are used for fitting the measured data. We increment the length of the time interval by 1~s within first 5~s, and then use the increment of 5~s up to 100~s. 
This allowed us to better examine the time-domain response of the considered circuits (discharge of both ladders),
connected as in Fig.~\ref{schemat_dl},
near the starting point $t=0$, and also their time-domain responses in long run,
which was in our case the interval up to 100 seconds. 
Discharges of the 60-steps domino ladder, 130-steps domino ladder, and nested domino ladder
are depicted in Figs.~\ref{fig:DL60relax}, \ref{fig:DL130relax}, and \ref{fig:NLrelax}, respectively.
The sampling period was $Ts=0.01$ s for all measurements of the discharges
used for the computations.

\begin{figure}[htbp]
\begin{center}
\includegraphics[width=\columnwidth]{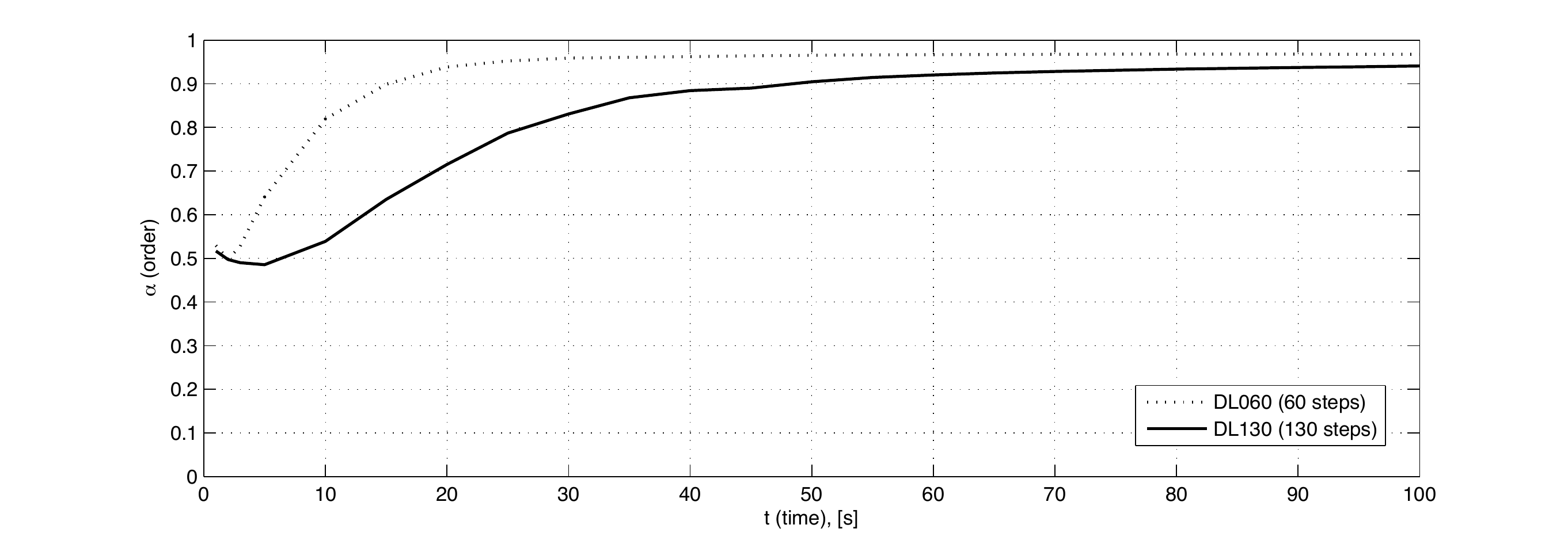}
\caption{Variable order $\alpha(t)$ for the 60-steps domino ladder (DL060), dotted line, and the 130-steps domino ladder (DL130), solid line.}
\label{fig:VO-plot}
\end{center}
\end{figure}

\begin{figure}[htbp]
\begin{center}
\includegraphics[width=\columnwidth]{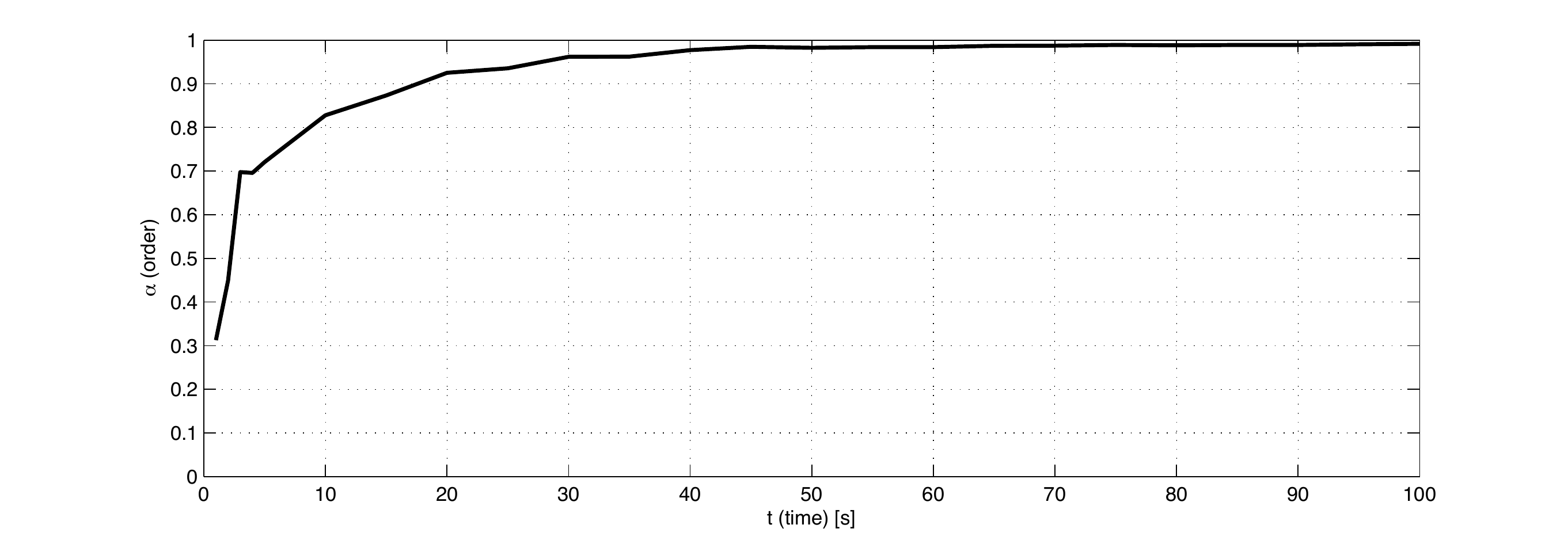}
\caption{Variable order $\alpha(t)$ for the nested domino ladder (NL14x14).}
\label{fig:VO-NestedLadder-plot}
\end{center}
\end{figure}

The results of these computations are presented in Table~\ref{tab:VO-table}
and in the Figs.~\ref{fig:VO-plot} and \ref{fig:VO-NestedLadder-plot}. 

The method of data fitting using the Mittag-Leffler function is implemented as a  Matlab routine~\cite{Podlubny-MLFFIT-FEX}, and the Mittag-Leffler function is computed also using our Matlab routine~\cite{Podlubny-MLF-FEX}.

\begin{table}[htb]
\centering
\caption{Variable order $\alpha(t)$ for the 60-steps domino ladder (DL060), 
the 130-steps domino ladder (DL130), and the nested ladder (NL14x14)} \label{tab:VO-table}
\begin{tabular}{| r | c | c| c|}
\hline
$t$ [s] & \multicolumn{3}{|c|}{ $\alpha(t)$} \\
\cline{2-4}
& \quad DL060 \quad & \quad  DL130 \quad  & \quad   NL14x14 \quad \\
\hline
1 & 0.5294 & 0.5167 & 0.3126 \\
2 & 0.4984 & 0.4972 & 0.4498 \\
3 & 0.5277 & 0.4901 & 0.6978 \\
4 &  0.5746 & 0.4821 & 0.6959 \\
5 & 0.6408 & 0.4855 &  0.7205\\
10 & 0.8195 & 0.5390 & 0.8278 \\
15 & 0.8986 & 0.6326 & 0.8732 \\
20 & 0.9385 & 0.7098 & 0.9249 \\
25 & 0.9523 & 0.7801 & 0.9354 \\
30 & 0.9586 & 0.8227 & 0.9618 \\
35 & 0.9604 & 0.8531 & 0.9620 \\
40 & 0.9620 & 0.8737 & 0.9770 \\
45 & 0.9638 & 0.8900 & 0.9847 \\
50 & 0.9651 & 0.9046 & 0.9823 \\
55 & 0.9661 & 0.9142 & 0.9837 \\
60 & 0.9668 & 0.9201 & 0.9837 \\
65 & 0.9670 & 0.9246 & 0.9868 \\
70 & 0.9674 & 0.9279 & 0.9872 \\
75 & 0.9678 & 0.9307 & 0.9886 \\
80 & 0.9677 & 0.9336 & 0.9879 \\
85 & 0.9676 & 0.9354 & 0.9888 \\
90 & 0.9677 & 0.9373 & 0.9887 \\
95 & 0.9672 & 0.9388 & 0.9904 \\
100 & 0.9672 & 0.9406 & 0.9915 \\
\hline
\end{tabular}
\end{table}

\section{Discussion}
\label{sec:discussion}

Our main conclusion is that both the domino ladder and the nested ladder
exhibit dual behavior in the frequency domain and in the time domain. 
In some frequency range or in some time interval they behave as 
fractional-order integrators of (almost) constant order. 
Outside of that frequency range or outside of that time interval
they behave as variable-order integrators; in one case that variable order
depends on frequency, in the other case it depends on the time.

The domino ladder behaves as Caputo-Weyl integrator of constant order $\alpha=0.5$ in a certain frequency range in the frequency domain.  
This means that in that frequency range it simply shifts the phase by $\alpha \pi / 2 = \pi/4$.
The frequency range where this behavior is observed can be made larger
by increasing the number of steps in the domino ladder. Outside of this frequency range
the domino ladder behaves as a variable-order system, where the order
depends on the frequency, as one can conclude directly from the Bode plots.

At the same time, in the time domain the same domino ladder behaves 
as an integrator of variable non-integer order, where the order depends on
the length of the time interval.  Close to the starting time instance $t=0$, 
the domino ladder behaves as an integrator of order $\alpha \approx 0.5$,
and with growing $t$ the domino ladder behaves closer and closer to
the classical integrator of order $1$. It should be mentioned that although 
the domino ladder order, $\alpha(t)$, tends to $1$, 
 the order $1$ is never reached. However, in many practical applications
 it is sufficient to neglect the transient effects for some initial time interval 
 and to assume that $\alpha(t) = 1$ for all $t$.

The similar observations hold for the nested ladder circuit, which has been
introduced in this paper. In some frequency range it behaves as an integrator
of order 0.25, and outside of this frequency range it behaves as 
a variable-order system, where the order depends on the frequency. 

In the time domain the nested ladder circuit behaves like a variable-order integrator, 
with order $\alpha(t)$ starting close to 0.25, and then increasing towards 1; 
the order $1$ is also never reached in the considered time interval. 

The frequency range and the time interval, where the order of the nested ladder 
is close to 0.25, can be extended by increasing the number of levels of the ladders
in the nested structure, and by increasing the numbers of steps in those ladders.

\section{Conclusions}
\label{sec:conclusions}

In this paper we have presented the experimental study 
of the two types of electrical circuits
made only of passive elements, which exhibit non-integer order behavior. 
One of them is the domino ladder, which already appeared in the works 
of other authors on the fractional-order systems. The other one is the circuit that we call 
the nested ladder and which was introduced in this paper.

For both these types of circuits we demonstrated that they 
should be considered not just as non-integer order systems, 
but as variable-order systems, where the order depends either on 
the frequency (in the frequency domain) or on the time variable
(in the time domain).  

While in the frequency domain the frequency-dependent variable order
is obvious directly from the Bode plots, providing the evidence of the 
variable-order behavior of the considered circuits in the time domain
required some additional tools. Namely, we suggested a method
of data fitting with the help of the Mittag-Leffler function, explained
a link between such fitting and fractional-order differential equations, 
and provided the Matlab routines for such fitting. 

The approach to identification of variable-order systems, that we presented 
in this paper, can be used for creating variable-order models for 
many other processes.

\section{Acknowledgment}

This entire work, including measurements and computations, 
has been done during the sabbatical stay of Dominik Sierociuk 
at the Institute of Control and Informatization of Processes,
BERG Faculty, Technical University of Kosice.


\begin{thebibliography}{99}


\bibitem{Podlubny}
I.~Podlubny,
\textit{Fractional Differential Equations},
San Diego: Academic Press,  1999.

\bibitem{Oldham}
K.~B.~Oldham and J.~Spanier,
\textit{The Fractional Calculus},
New York: Academic Press,  1974.


\bibitem{Magin-book}
R.~L.~Magin, 
\textit{Fractional Calculus in Bioengineering}, Connecticut: Begell House Publishers, 2006.


\bibitem{blas2010}
C. A.~Monje, Y. Q.~Chen, B. M.~Vinagre, D.~Xue and V.~Feliu,
 \emph{Fractional-order Systems and Controls}, 
 Series: Advances in Industrial Control, Springer, 2010.

\bibitem{Caponetto}
R.~Caponetto, G.~Dongola, L.~Fortuna, and I.~Petr\'a\v{s},
\textit{Fractional Order Systems:
Modeling and Control Applications}, Singapore: World Scientific, 2010.



\bibitem{petras2002}
I.~Petr\'a\v{s}, I.~Podlubny, P.~O'Leary, L. Dor{\v c}\'ak and B.~M. Vinagre, \emph{Analogue Realization of Fractional-Order Controllers},  TU Ko{\v s}ice: BERG Faculty,  Slovakia, 2002.

\bibitem{Petras2010mem}
I. Petr\'{a}\v{s}, 
 ``Fractional-order memristor-based Chua's circuit,"
\emph{IEEE Transactions on Circuits and Systems-II: Express Briefs}, vol. 57, no. 12,  pp.~975--979, 2010.

\bibitem{Petras2011}
 I.~Petr\'a\v{s},
\emph{Fractional-Order Nonlinear Systems: Modeling, Analysis and Simulation},
London: Springer and Beijing: HEP, 2011.

\bibitem{Westerlund2}
        S.~Westerlund and L.~Ekstam,
        ``Capacitor theory,"
        \emph{IEEE Trans. Dielectr. Electr. Insul.},
        vol.~1, no.~5, pp.~826--839, 1994.

\bibitem{PetrasETFA2009}
I. Petr\'{a}\v{s}, Y. Q. Chen and C. Coopmans,
``Fractional-order memristive systems," in
\emph{Proc. of the 14th IEEE International Conference ETFA'2009}, September 22-26, 2009, Mallorca, Spain. 

\bibitem{Schafer-k2}
I.~Schafer and K.~Kruger,
``Modelling of lossy coils using fractional derivatives,"
\textit{J. Phys. D: Appl. Phys.}, vol.~41, pp.~1--8, 2008.

\bibitem{Coopmans-Petras-Chen-k2}
	C.~Coopmans, I.~Petr\'a\v{s} and Y.~Q.~Chen,
	``Analogue Fractional-Order Generalized Memristive Devices," in
	\textit{Proc. of the ASME 2009 Conference},
              \#86861, San Diego, USA, August 30--September 2, 2009.
              
\bibitem{Dzelinski}
A. Dzielinski and D. Sierociuk, 
``Ultracapacitor modelling and control using discrete fractional order state-space model," 
\emph{Acta Montanistica Slovaca}, vol. 13, no. 1, pp. 136--145, 2008.

\bibitem{Dzelinski2}
A. Dzielinski, D. Sierociuk and G. Sarwas, 
``Some applications of fractional order calculus,"
\emph{Bulletin of The Polish Academy of Sciences-Technical Sciences},
 vol. 58, no. 4, pp. 583--592,  2010.         
              
\bibitem{Bohannan}
	G.~Bohannan, ``Analog Realization of a Fractional Control Element - Revisited,"  in
	\textit{Proc. of the IEEE Int. Conf. on Decision and Control}, 2002, Las Vegas, USA.
	
	
\bibitem{Sheng2011}
H. Sheng,  H. G. Sun, C. Coopmans, Y. Q. Chen  and G. W. Bohannan,
``A Physical experimental study of variable-order
fractional integrator and differentiator,"
\emph{Eur. Phys. J. Special Topics}, vol. 193,  pp. 93--104, 2011.




\bibitem{Nakagawa}
        M. Nakagawa and K. Sorimachi,
        ``Basic characteristics of a fractance device,"
        \emph{IEICE Trans. Fundamentals},
        vol.~E75 - A, no.~12, pp.~1814--1818, 1992.
   
 
\bibitem{Carlson2}
        G. E. Carlson and C. A. Halijak,
        ``Approximation of fractional capacitors $(1/s)^{1/n}$
        by a regular Newton process,"
        \emph{IEEE Trans. on Circuit Theory}, vol.~11, no.~2,  pp.~210--213, 1964.

        

\bibitem{Petras-Analog}
I. Podlubny,  I. Petr\'{a}\v{s},  B. Vinagre, P. O'Leary and  L. Dor{\v c}\'ak, 
``Analogue Realizations of Fractional-order Controllers", 
\emph{Nonlinear Dynamics}, vol. 29, no.~1--4, pp. 281--296, 2002.

\bibitem{erfani2002}
S.~Erfani, ``Evaluation and realization of continued-fraction expansion revisited," 
\emph{Computers and Electrical Engineering}, vol.  28, pp. 311--316, 2002.

\bibitem{Biswas}
K. Biswas, S. Sen and P. K. Dutta,
``Realization of a Constant Phase Element and Its Performance Study in a Differentiator Circuits,"
\emph{IEEE Trans. on Circuits and Systems II: Express Brief}, vol. 53, no. 9, pp. 802--806, 2006. 

\bibitem{roy1967} 
S. D.~Roy, ``On the Realization of a Constant-Argument Immittance or Fractional Operator," 
\emph{ IEEE Transactions on Circuit Theory}, vol.  14, no.~3, pp. 264--274, September 1967.


\bibitem{roy1974} 
S. D.~Roy, ``Constant argument immittance realization by a distributed RC network," 
\emph{IEEE Transactions on Circuits and Systems},  vol.  21, no. 5, pp. 655--658, September 1974.

\bibitem{Krishna2008} 
B. T.~Krishna and K. V. V. S.~Reddy, ``Active and Passive Realization of Fractance Device of Order 1/2," 
\emph{Active and Passive Electronic Components}, vol. 2008, article ID 369421, 5 pages.


\bibitem{Yifei2005} 
P.~Yifei, Y.~Xiao, L.~Ke, Z.~Jiliu, Z~Ni, Z.~Yi and P.~Xiaoxian,
``Structuring analog fractance circuit for 1/2 order fractional calculus," in 
\emph{Proc. of The 6th International Conference on ASICON 2005}, vol. 2, 24--30 October, 2005, pp. 1136--1139. 

\bibitem{Wang}
        J. C. Wang,
        ``Realizations of generalized Warburg impedance with
        RC ladder networks and transmission lines,"
        \emph{J. of Electrochem. Soc.}, vol.~134, no.~8,  pp.~1915--1920, August 1987.

\bibitem{Manabe}
        S.~Manabe,
       ``The Non-Integer Integral and its Application to Control Systems,"
	  \emph{ETJ of Japan}, vol. 6, no. 3-4, pp.~83--87, 1961.

\bibitem{Oldham73}
K.~B.~Oldham, 
``Semiintegral Electroanalysis: Analog Implementation,"
\emph{Analytical Chemistry}, vol. 45, no. 1, pp. 39--47, 1973.


\bibitem{SierociukMMAR11}
D. Sierociuk and A. Dzielinski,
``New method of fractional order integrator analog
modeling for orders $\alpha = 0.25$ and $\alpha= 0.5$," in
\emph{Proc. of the 16th International Conference on
Methods and Models in Automation and Robotics},
22 -- 25 August 2011, Amber Baltic Hotel, Miedzyzdroje, Poland (to be published).

\bibitem{Podlubny-MLFFIT-FEX}
I. Podlubny. \textit{Fitting experimental data using the Mittag-Leffler function.}
Matlab Central File Exchange, submission  \#32170.
Available at 
\url{http://www.mathworks.com/matlabcentral/fileexchange/32170}.


\bibitem{Podlubny-MLF-FEX}
I. Podlubny and M. Kacenak. \textit{Mittag-Leffler function.} 
Matlab Central File Exchange, submission \#8738. 
Available at 
\url{http://www.mathworks.com/matlabcentral/fileexchange/8738}.


\end{thebibliography}
\end{document}